\documentclass[letterpaper, 10 pt, conference]{ieeeconf}  

\usepackage{amsmath,amssymb,amsfonts}
\usepackage{graphicx}
\usepackage{etoolbox}
\usepackage{cases}

\usepackage{algorithmic}
\usepackage{graphicx}
\usepackage{textcomp}
\usepackage{xcolor}
\usepackage{array}

\usepackage{mathrsfs}
\newtheorem{theorem}{Theorem}
\newtheorem{lemma}{Lemma}[section]
\newtheorem{remark}{Remark}
\newtheorem{definition}{Definition}
\newtheorem{property}{Property}

\usepackage{comment}
\usepackage{tabularx}
\usepackage{multirow}
\IEEEoverridecommandlockouts                              

\title{\LARGE \bf
An online optimization algorithm for tracking a linearly varying optimal point with zero steady-state error}

\author{Alex (Xinting) Wu, Ian R. Petersen, Valery Ugrinovskii, Iman Shames
\thanks{*This work was supported by the Australian Research Council under grants DP200102945, DP210102454 and DP230102443.
}
\thanks{$^{1}$A. Wu, I. R. Petersen and I. Shames are with the CIICADA Lab, School of Engineering, The Australian National University, Canberra, ACT 2601,
Australia (e-mail: u5847417@anu.edu.au; i.r.petersen@gmail.com; iman.shames@anu.edu.au).}%
\thanks{$^{2}$V. Ugrinovskii was with the School of Engineering and IT, UNSW Canberra, Canberra, ACT 2600, Australia.}%
}

\begin{document}

\maketitle
\thispagestyle{empty}
\pagestyle{empty}

\begin{abstract}
In this paper, we develop an online optimization algorithm for solving a class of nonconvex optimization problems with a linearly varying optimal point. The global convergence of the algorithm is guaranteed using the circle criterion for the class of functions whose gradient is bounded within a sector. Also, we show that the corresponding Luré-type nonlinear system involves a double integrator, which demonstrates its ability to track a linearly varying optimal point with zero steady-state error. The algorithm is applied to solving a time-of-arrival based localization problem with constant velocity and the results show that the algorithm is able to estimate the source location with zero steady-state error.

\end{abstract}

\section{INTRODUCTION}
Online optimization \cite{Introonline} has emerged as a pivotal area of research with significant implications across various fields, including real-time control systems, signal processing, and machine learning. An online optimization problem typically involves making a sequence of decisions in real-time, where the objective function changes over time. The primary challenge lies in solving an optimization problem with an objective function that shifts as time progresses, while ensuring that the time-varying solution converges with minimal steady-state error, thereby maintaining high accuracy and stability over time. Recently, the design of algorithms to solve such problems has been investigated in \cite{Zampieri2023}, employing tools from robust control theory, especially in solving convex problems. The paper \cite{Zampieri2023} designed control-based online optimization algorithms based on the internal model principle, which states that for a control system to eliminate steady-state error, the controller must incorporate a model of corresponding disturbances and input signals \cite{FRANCIS1976457}. This improves the accuracy and performance of the solutions. In this paper, we consider a special case of the problems considered in \cite{Zampieri2023} in which the optimal point of the objective function varies linearly with time. Then, as opposed to the approach of \cite{Zampieri2023}, we design an algorithm which minimizes the root convergence rate \cite{IterativeSolution} along with ensuring zero steady-state error; see also  \cite{AlexACC2024,9745160,UGRINOVSKII2023111129}. This leads to a simple closed form formula for the resulting algorithm.

We consider an unconstrained optimization problem:
\begin{equation}
\label{mini}
   x^*(t) = \arg \min f(x,t),
\end{equation}
where $f: \mathbb{R}^{n+1\times n+1} \rightarrow \mathbb{R}$ is the cost function that attains its minimum at $x^*(t)$. In the case which $f(x,t)$ is independent of $t$, \cite{AlexACC2024} proposes an approach to this problem through an algorithm which we generalize to the following form:
\begin{align}
    \label{uprule}
    x(t+1) &= x(t) -\alpha \nabla_x f(x(t),t) + \gamma \nabla_x f(x(t-1),t) \nonumber \\& \quad+ \beta(x(t)-x(t-1)),
\end{align}
where $x(t) \in\mathbb{R}^n$ is the current estimate, $t$ is the iteration index, $\alpha$ is the step size, $\beta$ is a momentum parameter, $\gamma$ is an additional parameter, and $\nabla_x f(x,t)$ is the gradient of the cost function.

In order to obtain zero steady-state error with a linearly varying optimal point, we introduce an extra integration into the algorithm. This implies we restrict our attention to the special case of (\ref{uprule}) where $\beta=1$. Then the recursion (\ref{uprule}) becomes
\begin{align}
    \label{uprulebeta1}
    x(t+1) &= x(t) -\alpha \nabla_x f(x(t),t) + \gamma \nabla_x f(x(t-1),t)\nonumber \\& \quad + x(t)-x(t-1).
\end{align}
Also, we consider the class of possibly nonconvex cost functions as defined in \cite{AlexACC2024}.

\begin{definition}
    Given $L \geq m >0$, let $\mathscr{F}_{m,L}$ denote the class of cost functions $f(x,t)$ with the following property:
    
    The function $f(x,t)$ is twice differentiable such that $f(x,t) = \Tilde{f} \left(x(t)-x^*(0)-at \right)$ and for every $x(t) \in \mathbb{R}^n$,
\begin{align}
\label{classfunc}
    &\left( m \left(x(t)-x^*(t) \right)- \nabla_x f(x(t),t) \right)^T \nonumber \\& \times \left( L \left(x(t) \nonumber -x^*(t) \right)-\nabla_x f(x(t),t) \right)\nonumber \\ & \leq 0, \quad \text{for all} \quad t\geq0.
\end{align}

Here, $x^*(t) = x^*(0)+at$ is the linearly varying optimal point of $f(x,t)$, where $a$ is a constant vector. 
\end{definition}
Note that for twice differentiable functions, (\ref{classfunc}) implies that
\begin{equation}
\label{classfuncmL}
    mI \leq \nabla_x^2 f(x^*(t),t) \leq LI,\quad \text{for all} \quad t\geq0.
\end{equation}

The inequality (\ref{classfunc}) is a multivariate analogue of a sector bound on the gradient of a scalar function, where the gradient $\nabla_x f(x,t)$ belongs to the sector $[m,L]$ \cite{haddad2008nonlinear}. 

In this paper, we develop an online optimization algorithm of the form (\ref{uprulebeta1}) to solve a class of optimization problems with cost functions $f(x,t) \in \mathscr{F}_{m,L}$. We write the algorithm (\ref{uprulebeta1}) as a Luré-type nonlinear system where the corresponding transfer function incorporates a double integrator. The system is recognized as a Type 2 system, in which the system will have zero steady-state error in response to a ramp input \cite[Table 4-4, p.200]{10.5555/201033}. Moreover, the absolute stability of the corresponding system is established by employing the circle criterion, which guarantees the global convergence of the algorithm. We specifically consider using the root convergence factor \cite[p. 288]{IterativeSolution} as a measure of the rate of convergence. We also apply the algorithm to a time-of-arrival (TOA) based localization problem, demonstrating that the algorithm is able to estimate the linearly varying source location with zero steady-state error.

This paper is organized as follows: Section \ref{sec2} presents background material and preliminary results. In Section \ref{sec3}, we describe and prove our main result. Section \ref{sec4} illustrates our algorithm with an example of a TOA based localization problem. Finally, the paper is concluded in Section \ref{sec5}.

\section{Background and Preliminary Results} \label{sec2}

\subsection{Heavy Ball Method}
One of the most well known methods to solve the optimization problem (\ref{mini}) is the heavy ball method \cite{POLYAK19641}. This method enhances the standard gradient descent method by incorporating a momentum term leading to an optimization algorithm of the form
\begin{equation}
\label{HB}
    x(t+1) = x(t) -\alpha \nabla f(x(t),t) + \beta(x(t)-x(t-1)).
\end{equation}
The heavy ball method combines the current gradient with information of the previous step to accelerate the convergence. With Polyak's choice of parameters, the method demonstrates fast local convergence. This choice of parameters is given by 
\begin{equation}
\label{polypara}
    \begin{aligned}
        \alpha = \alpha_{\text{\scriptsize HB}} =&\frac{4}{(\sqrt{L}+\sqrt{m})^2}= \frac{4}{m(1+\sqrt{\kappa})^2},\\
        \beta = \beta_{\text{\scriptsize HB}} =& \frac{(\sqrt{L}-\sqrt{m})^2}{(\sqrt{L}+\sqrt{m})^2}= \frac{(\sqrt{\kappa}-1)^2}{(\sqrt{\kappa}+1)^2}.
    \end{aligned}
\end{equation}
The global convergence property of the heavy ball method has been studied in \cite{9745160} for the class of cost functions $\mathscr{F}_{m,L}$. The paper \cite{9745160} has shown that for any $\kappa = L{/}m<\kappa_{\text{\scriptsize HB}} = 3+2\sqrt{2}$, the heavy ball method achieves global asymptotic convergence for any function $f(x,t)$ $\in \mathscr{F}_{m,L}$, when $f(x,t)$ is independent of $t$.

\subsection{Generalized accelerated gradient method}
The generalized accelerated gradient method of \cite{AlexACC2024} extends the method of \cite{9745160} by introducing an additional parameter, leading to an algorithm of the form (\ref{uprule}). The corresponding choice of parameters $\alpha$, $\beta$, $\gamma$ in (\ref{uprule}) are denoted by  $\alpha = \alpha_{\text{\scriptsize GAG}}$, $\beta = \beta_{\text{\scriptsize GAG}}$ and $\gamma = \gamma_{\text{\scriptsize GAG}}$. The algorithm presented in \cite{AlexACC2024} exhibits a globally convergent property, and the introduction of this additional parameter yields an enhanced R-convergence rate for any cost function $f(x,t)$ $\in \mathscr{F}_{m,L}$, when $f(x,t)$ is independent of $t$.

\subsection{Preliminary Results}
Let $\Sigma(\alpha,\gamma,m,L)$ denote the algorithm corresponding to the recursion (\ref{uprulebeta1}) depending on $\alpha, \gamma$ and the set of cost functions $\mathscr{F}_{m,L}$ defined by (\ref{classfunc}) where $L>m>0$.

\begin{definition}(\cite[p. 765]{haddad2008nonlinear})
    Given $\alpha$, $\gamma$, $m$ and $L$, the algorithm $\Sigma(\alpha,\gamma,m,L)$ is said to be globally asymptotically convergent if given any $f(x,t) \in \mathscr{F}_{m,L}$ with optimal points $x^*(t)$, the following conditions hold:
    \begin{enumerate}
        \item $x(t)=x^*(t)$ is a Lyapunov stable solution of (\ref{uprule}). That is, for all $\sigma>0$ there exists $\delta=\delta(\sigma)>0$ such that if $\|x(0)-x^*(0)\|<\delta$ and $\|x(1)-x^*(1)\|<\delta$, then $\|x(t)-x^*(t)\|<\sigma$, $t = 2, 3, ...$.\\
        \item The algorithm is globally asymptotically convergent. That is, for any $x(0), x(1)\in \mathbb{R}^n$, $\lim_{t \rightarrow \infty}(x(t)-x^*(t))=0$.
    \end{enumerate}
\end{definition}

\begin{remark}
Note that if an algorithm $\Sigma(\alpha,\gamma,m,L)$ is globally asymptotically convergent according to this definition, then it solves the optimization problem (\ref{mini}) with zero steady-state error since $\lim_{t \rightarrow \infty}(x(t)-x^*(t))=0$.
\end{remark}

The following definition describes the root convergence factor, also referred to as the R-convergence rate.
\begin{definition}(\cite[p. 288]{IterativeSolution})
    Given $\alpha$, $\gamma$, $m$ and $L$ such that $\Sigma(\alpha,\gamma,m,L)$ is globally asymptotically convergent, then the R-convergence rate is defined by
    \begin{equation}
    \label{rrate}
        r(\alpha,\gamma,m,L) = \sup_{x(0),x(1),f(x,t) \in \mathscr{F}_{m,L}} \lim_{t \rightarrow \infty} \|x(t)-x^*(t)\|^{\frac{1}{t}},\nonumber
    \end{equation}
    where $x(t)$ is the solution to (\ref{uprule}) corresponding to $f(x,t) \in \mathscr{F}_{m,L}$ with optimal point $x^*(t)$ and initial conditions $x(0)$, $x(1)$.
\end{definition}

Our analysis leverages absolute stability theory, thereby formulating the algorithm (\ref{uprulebeta1}) as a Luré type system. To achieve this, first define
\begin{align}
    w(t) & \triangleq  x(t-1) + \gamma u(t-1) \label{zt},\\
    u(t) & \triangleq -\nabla f(y,t) \label{ut},
\end{align}
where
\begin{equation}
    y(t) = x(t).\label{ytxt}
\end{equation}
From (\ref{zt}) and (\ref{uprule}), it follows that
\begin{equation}
    w(t+1) = x(t) +\gamma u(t) \label{wt+1}
\end{equation}
and
\begin{equation}
    x(t+1) =-w(t)+2x(t) +\alpha u(t) \label{yt+1}.
\end{equation}
Equations (\ref{ut}), (\ref{ytxt}), (\ref{wt+1}) and (\ref{yt+1}) can be written as a Luré system with a standard state space form as
\begin{align}
\label{ssm2}
    X(t+1) =& AX(t) + B u(t), \nonumber\\
    y(t) =& C x(t), \\
    u(t) =& -\nabla f(y,t), \nonumber
\end{align}
where $X(t) = [w(t)^T \quad  x(t)^T]^T $,
\begin{equation}
\label{itermatr}
     A = A_0 \otimes I_n, \quad B = B_0 \otimes I_n, \quad  C = C_0 \otimes I_n,
\end{equation}
and
\begin{equation}
\label{linmatr}
    A_0 = \begin{bmatrix}
        0&1 \\ -1 & 2
    \end{bmatrix},
    B_0 = \begin{bmatrix}
        \gamma \\ \alpha
    \end{bmatrix},
    C_0 = \begin{bmatrix}
        0&1
    \end{bmatrix}.
\end{equation}
Here $\otimes$ denotes the Kronecker product.

Using the the circle criterion \cite{haddad2008nonlinear}, the absolute stability of the Luré system (\ref{ssm2}), (\ref{itermatr}), (\ref{linmatr}) can be secured. Let $G(z)$ denote the transfer function of the linear part of the system (\ref{ssm2}),  (\ref{itermatr}), (\ref{linmatr}). Then 
\begin{equation*}
    G(z) = G_0(z) I_n
\end{equation*}
where
\begin{align}
    G_0(z) &=  C_0(zI-A_0)^{-1}B_0 \nonumber \\
    &=  \frac{\alpha z  -\gamma}{(z-1)^2}.
\end{align}
The following figure shows the block diagram of the system  (\ref{ssm2}), (\ref{itermatr}), (\ref{linmatr}).
\begin{figure}[!htb] 
    \centering
    \includegraphics[scale=0.43]{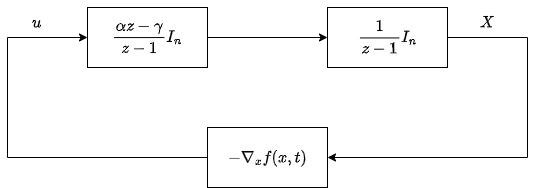}
    \caption{Block diagram corresponding to  Luré system representation}
    \label{fig:BD}
    \vspace{-0.2cm}
\end{figure}

In Figure \ref{fig:BD}, the block diagram depicts the presence of a double integrator. This component signifies the system's capability to effectively track a linearly varying optimal point.  

Define $\Tilde{A}(\Delta)=A-B \Delta C$, where $A, B, C$ are defined in (\ref{itermatr}), and $\Delta$ is a symmetric matrix which satisfies (\ref{classfuncmL}). We also define $\Tilde{A_0}(\lambda) = A_0- \lambda B_0 C_0$, where $m \leq \lambda \leq L$.

The following lemma relates the R-convergence rate and the spectral radius of $\Tilde{A_0}(\lambda)$.
\begin{lemma}
\label{rconverge}
    Given any $\alpha$, $\gamma$, $m$, $L$, such that the algorithm $\Sigma(\alpha,\gamma,m,L)$ is globally asymptotically convergent,
    \begin{align}
        r(\alpha,\gamma,m,L) =&\sup_{mI \leq \Delta \leq LI} \rho(\Tilde{A}(\Delta)) \label{matrixr}\\
        =&\sup_{m \leq \lambda \leq L} \rho(\Tilde{A_0}(\lambda)). \label{scalarr}
    \end{align}
    where $\rho(\cdot)$ denote the spectral radius of its matrix argument.
\end{lemma}
\begin{proof}
    Let 
    \begin{align} 
    \label{xttil}
        \Tilde{x}(t) &= x(t)-x^*(t) \nonumber \\
        &= x(t) - x^*(0)-at
    \end{align}
    Then
    \begin{equation}
    \label{ftoftil}
        f(x,t) = \Tilde{f} \left(x(t)-x^*(0)-at \right) =  \Tilde{f}(\Tilde{x}(t)).
    \end{equation}
    Also,
    \begin{equation}
    \label{tilnabla}
        \nabla_x f(x,t) = \nabla_x \Tilde{f}\left(x(t)-x^*(0)-at \right) = \nabla_x \Tilde{f}(\Tilde{x}(t)).
    \end{equation}
    It follows from (\ref{xttil})
    \begin{equation}
    \label{xt+1til}
        \Tilde{x}(t+1) = x(t+1)-x^*(0)-a(t+1).
    \end{equation}
    Substituting (\ref{uprulebeta1}) into (\ref{xt+1til}) yields
    \begin{align}
        \Tilde{x}(t+1) &=x(t) -\alpha \nabla_x f(x(t),t) + \gamma \nabla_x f(x(t-1),t) \nonumber \\& \quad + x(t)-x(t-1)-x^*(0)-a(t+1).
    \end{align}
    It follows from (\ref{xttil}) and (\ref{tilnabla}) that
    \begin{align}
    \label{upruletilx1}
        \Tilde{x}(t+1) &= \Tilde{x}(t) +x^*(0)+at -\alpha \nabla_x \Tilde{f}(\Tilde{x}(t)) \nonumber \\ & \quad + \gamma \nabla_x \Tilde{f}(\Tilde{x}(t-1)) + \Tilde{x}(t) +x^*(0) \nonumber \\ & \quad +at  -(\Tilde{x}(t-1) +x^*(0)+a(t-1)) \nonumber \\ & \quad -x^*(0)-a(t+1).
    \end{align}
    Simplifying (\ref{upruletilx1}) yields
    \begin{align}
    \label{upruletilx}
       \Tilde{x}(t+1) &= \Tilde{x}(t)-\alpha \nabla_x \Tilde{f}(\Tilde{x}(t)) + \gamma \nabla_x \Tilde{f}(\Tilde{x}(t-1)) \nonumber\\& \quad + \Tilde{x}(t)-\Tilde{x}(t-1).
    \end{align}
    In this way, we are able to write the algorithm $\Sigma(\alpha,\gamma,m,L)$ in terms of recursion (\ref{upruletilx}) and its the corresponding Luré system system (\ref{ssm2}), (\ref{itermatr}), (\ref{linmatr}) is transformed into a time-invariant system related to the recursion (\ref{upruletilx}). Therefore, the conditions (\ref{matrixr}), (\ref{scalarr}) follow from \cite[Theorem 10.1.4, p. 301]{IterativeSolution} and \cite[Theorem 1]{UGRINOVSKII2023111129}. 
\end{proof}

In order to apply the formula (\ref{scalarr}) for the R-convergence rate, we consider the characteristic polynomial of the matrix  $\Tilde{A_0}(\lambda)$. This characteristic polynomial is given by 
\begin{align}
\label{charpoly}
    \chi (z) =& \det (zI-\Tilde{A_0}(\lambda)) \nonumber \\
    =&\det \begin{bmatrix}
        z & \lambda \gamma-1 \\
        1 & z+\lambda \alpha-2 \nonumber
    \end{bmatrix}\\
    =& z^2+(-2 +\lambda \alpha)z-\lambda \gamma+1 \nonumber\\
    =&z^2+a_1 z+ a_2,
\end{align}
where
\begin{equation*}
    a_1 = -2+\lambda \alpha, \quad
    a_2 = 1 - \lambda \gamma.
\end{equation*}

The following result is a version of the circle criterion applied to the system (\ref{ssm2}), (\ref{itermatr}), (\ref{linmatr}).
\begin{lemma}
\label{circlecr}
    For given $\alpha$, $\gamma$, $m$, $L$, the algorithm $\Sigma(\alpha,\gamma,m,L)$ is globally asymptotically convergent if the following conditions are satisfied:
    \begin{enumerate}
        \item $1+mG_0(z) \neq 0$ for $|z| \geq 1$;
        \item the transfer function $H_0(z)$ 
            \begin{align}
                H_0(z) &= [1+LG_0(z)][1+mG_0(z)]^{-1} \nonumber\\
                &= \frac{z^2+z( L\alpha  -2)- L\gamma+1}{z^2+z(m\alpha -2)-m\gamma + 1} \label{H0zml}
            \end{align}
        is strictly positive real \cite[Definition 5.18]{haddad2008nonlinear}.
    \end{enumerate}
\end{lemma}
\begin{proof}
    It follows from the proof of Lemma \ref{rconverge} that the system  (\ref{ssm2}), (\ref{itermatr}), (\ref{linmatr}) can be transformed into a time-invariant system corresponding to the recursion (\ref{upruletilx}). Then, using the properties of the Kronecker product, the result follows from the discrete-time circle criterion; e.g., see Problem 13.27 in \cite{haddad2008nonlinear}. 
\end{proof}

If the conditions of the above lemma are satisfied, the algorithm $\Sigma(\alpha,\gamma,m,L)$ is said to be \textbf{globally asymptotically convergent via the circle criterion}.

\section{Main Result} \label{sec3}
We aim to develop an algorithm of form (\ref{uprulebeta1}) which is globally asymptotically convergent via the circle criterion with R-convergence rate $0<\rho<1$.

\begin{theorem}
\label{maintheo}
Given any $m>0$ and $L>m>0$, let $\kappa = \frac{L}{m}$. Then there exists an algorithm $\Sigma(\alpha,\gamma,m,L)$ which is globally asymptotically convergent via the circle criterion with R-convergence rate 
\begin{align*}
    r(\alpha^*,\gamma^*,m,L)&=\rho^*(\kappa)\\
    &= \sqrt{\frac{\kappa-1}{\kappa+1}}.
\end{align*}
Also, the corresponding values of $\alpha$ and $\gamma$ are given by
\begin{equation*}
      \alpha^* = \frac{2}{m\kappa}, \quad \gamma^* = \frac{1}{m(\kappa+1)}.  
\end{equation*}

In addition, given any other algorithm of the form (\ref{uprulebeta1}), which is globally asymptotically convergent via the circle criterion, then its R-convergence rate satisfies 
\begin{equation*}
    r(\alpha,\gamma,m,L) \geq \rho^*(\kappa).
\end{equation*}
\end{theorem}

\begin{remark}
    Since $\kappa = \frac{L}{m}$
    \begin{equation*}
        r(\alpha^*,\gamma^*,m,L)=\frac{\sqrt{L-m}}{\sqrt{L+m}}
    \end{equation*}
    and the corresponding values of $\alpha$ and $\gamma$ are
    \begin{equation*}
      \alpha^* = \frac{2}{L}, \quad \gamma^* = \frac{1}{m+L}.  
\end{equation*}
\end{remark}
In order to prove Theorem \ref{maintheo}, we aim to construct an algorithm $\Sigma(\alpha,\gamma,m,L)$ such that the following properties are satisfied for suitable $0<\rho<1$:
\begin{property}
\label{prop1}
    The characteristic polynomial $\chi(z)$ defined in (\ref{charpoly}) has all its roots inside the circle of radius $0<\rho<1$ for all $m\leq \lambda \leq L$;
\end{property}
\begin{property}
\label{prop2}
    $1+mG_0(z) \neq 0$ for $|z| \geq 1$ and the transfer function $H_0(z)$ (\ref{H0zml}) is strictly positive real.
\end{property}

Note that the numerator of $1+mG_0(z)$ is identical to the characteristic polynomial $\chi(z)$ with $\lambda=m$. Hence, because $0<\rho<1$, Property \ref{prop1} will imply that the first part of Property \ref{prop2} is satisfied.

Property \ref{prop1} ensures the algorithm $\Sigma(\alpha,\gamma,m,L)$ has an R-convergence rate of $0<\rho<1$. This follows from Lemma \ref{rconverge}. Property \ref{prop2} guarantees that the algorithm is asymptotically convergent via circle criterion, according to Lemma \ref{circlecr}. 

\begin{figure}[!htb] 
\begin{center}
\vspace{-0.4cm}
\includegraphics[scale=0.65]{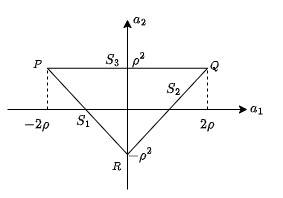}
\end{center} 
\vspace{-0.5cm}
\caption{Triangle shaped region $T_\rho$.}
\vspace{-0.5cm}
\label{fig:triangleTp}
\end{figure}

We first consider some results related to Property \ref{prop1}. These results parallel the corresponding results given in \cite{AlexACC2024} but are repeated here for the sake of completeness. We apply the Jury test \cite{10.5555/201033} to a scaled version of the characteristic polynomial $\chi(z)$ (\ref{charpoly}):
\begin{equation*}
\label{jury}
    \chi \left( \frac{z}{\rho} \right) = \left( \frac{z}{\rho} \right) ^2+a_1 \frac{z}{\rho}+ a_2, 
\end{equation*}
for $0<\rho<1$. This leads to a triangle shaped region in the  $(a_1,a_2)$-plane as shown in Figure \ref{fig:triangleTp}. This region is such that Property \ref{prop1} is satisfied if and only if the point $(a_1,a_2)$ lies inside it.
\\
Let $T_\rho$ denote the triangle shown in Figure \ref{fig:triangleTp}, 
\begin{equation*}
    T_\rho = \text{conv}\{\text{P,Q,R}\},
\end{equation*}
where $\text{conv} \{\cdot \}$ denotes the convex hull and $P = (-2\rho,\rho^2), Q = (2\rho,\rho^2)$ and $R = (0,-\rho^2)$.\\
The boundary of the region $T_\rho$ in Figure \ref{fig:triangleTp} is represented by
\begin{equation*}
\label{nablaT}
    \partial T_\rho = T_\rho \setminus \text{int}(T_\rho) = S_1 \cup S_2 \cup S_3,
\end{equation*}
where $S_1 = \text{conv} \{\text{P,R}\}, S_2= \text{conv} \{\text{Q,R}\}$ and $S_3 =\text{conv} \{\text{P,Q}\} $.
We define a line segment $J(\alpha,\gamma,m,L)$ with slope $-\frac{\gamma}{\alpha}$ inside $T_\rho$ as follows:
\begin{multline}
    \label{lineI}
    J(\alpha,\gamma,m,L) = \{(a_1,a_2): a_1 = -2+\lambda \alpha, \\a_2 = 1 - \lambda \gamma,\lambda \in [m,L]\}.
\end{multline}
An example of such a line segment $J(\alpha,\gamma,m,L)$ is shown in Figure \ref{fig:case1}. This line segment represents the set of characteristic polynomials $\chi(z)$ defined in (\ref{charpoly}) as the parameter $\lambda$ ranges over the interval $\lambda\in [m,L]$. We can write
\begin{equation*}
    J(\alpha,\gamma,m,L) = \text{conv}\{E_1(\alpha,\gamma,m),E_2(\alpha,\gamma,L)\},
\end{equation*}
where
\begin{align}
    E_1 (\alpha,\gamma,m) &= (-2+m\alpha,1-m\gamma), \label{E1}\\
    E_2 (\alpha,\gamma,L) &= (-2+L\alpha ,1-L\gamma) \label{E2}.
\end{align}
We now define the following quantities for $\gamma \geq 0$:
\begin{align}
    \Bar{\kappa}(\alpha,\gamma) &= \sup \{\kappa: J(\alpha,\gamma,m,L) \subset T_\rho \}, \label{ksupk}\\
    \Bar{m}(\alpha,\gamma) &= \inf\{m>0 : E_1 (\alpha,\gamma,m) \in T_\rho \} \label{barm},\\
    \Bar{L}(\alpha,\gamma) &= \sup \{ L>0:  E_2 (\alpha,\gamma,L) \in T_\rho \} \label{barL}.
\end{align}

The quantity $\Bar{\kappa}(\alpha,\gamma)$ represents the maximum value of the ratio $\kappa=\frac{L}{m}$ such that the interval $J(\alpha,\gamma,m,L)$ is contained in $T_\rho$. Using Lemma \ref{rconverge}, this will be the maximum value of $\kappa=\frac{L}{m}$ such that an algorithm $\Sigma(\alpha,\gamma,m,L)$ has R-convergence $\rho$. The quantities $\Bar{m}(\alpha,\gamma)$ and $\Bar{L}(\alpha,\gamma)$ will be used in the calculation of $\Bar{\kappa}(\alpha,\gamma)$. Indeed, it follows from the above definitions that 
\begin{equation*}
    \Bar{\kappa}(\alpha,\gamma) = \frac{\Bar{L}(\alpha,\gamma)}{\Bar{m}(\alpha,\gamma)}.
\end{equation*}

\begin{lemma}
\label{lem1}
Given $\alpha,\gamma$, suppose $\Bar{L}=\Bar{L}(\alpha,\gamma)>0$. Then $E_2(\alpha,\gamma,\Bar{L}) \in \partial T_\rho$.
\end{lemma}
\begin{proof}
Suppose $E_2(\alpha,\gamma,\Bar{L}) \not\in \partial T_\rho$. Then $E_2(\alpha,\gamma,\Bar{L}) \in \text{int}(T_\rho)$. Since $E_2(\alpha,\gamma,L)$ is continuous in $L$, there exists an $\hat{L}>\Bar{L}$ so that $E_2(\alpha,\gamma,\hat{L}) \in \text{int}(T_\rho)$. This contradicts the definition of $\Bar{L}$ in equation (\ref{barL}). 
\end{proof}

For a given value of $\gamma \geq 0$, we now consider the maximum achievable value of $\Bar{\kappa}(\alpha,\gamma)$. Considering the following definition:
\begin{equation}
\label{kbarmax}
    \bar \kappa_{\text{max}}(\gamma) = \sup_{\alpha} \Bar{\kappa}(\alpha,\gamma),
\end{equation}
$\bar \kappa_{\text{max}}(\gamma)$ represents the largest value of $\kappa$ for which there exists an algorithm $\Sigma(\alpha,\gamma,m,L)$ with the R-convergence rate $\rho$, when $\gamma$ is given and fixed. This is shown in the following lemma.

\begin{lemma}
\label{lemkmaxp1}
For a given $\gamma \geq 0$, the supremum on the right hand side of (\ref{kbarmax}) is attained and is equal to
\begin{equation}
\label{kmaxp}
    \bar \kappa_{\text{max}}(\rho) = \frac{1+\rho^2}{1-\rho^2}.
\end{equation}
The values of $\alpha, \gamma$ which attain the supremum are given by
\begin{align}
    \alpha (\rho) &= \frac{1}{m} \left(\frac{2(1-\rho^2)}{1+\rho^2} \right), \label{kmaxa} \\
    \gamma (\rho) &=\frac{1}{m}(1-\rho^2). \label{kmaxy}
\end{align}
\end{lemma}
\begin{proof}
In order to prove this lemma, we consider the following two cases.

\textbf{Case 1}. Let $E_1 = (-2\rho+\sigma_1,\rho^2-\sigma_2)$, $\sigma_1 \in [0,2\rho]$, $\sigma_2\in [0,\rho\sigma_1]$ and $E_2\in S_2 $ as shown in Figure \ref{fig:case1}.
\begin{figure}[!htb] 
\begin{center}
\vspace{-0.4cm}
\includegraphics[scale=0.65]{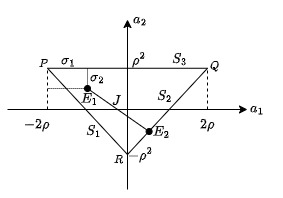}
\end{center} 
\vspace{-0.8cm}
\caption{Illustration of Case 1.}
\vspace{-0.2cm}
\label{fig:case1}
\end{figure}

It follows from (\ref{E1}) that $E_1 = (-2+m\alpha,1-m\gamma)$ and hence
\begin{align}
    -2\rho+\sigma_1 &= -2+m\alpha \label{c1e1x},\\
    \rho^2-\sigma_2 &=1-m\gamma. \label{c1e1y}
\end{align}
In addition, on the line segment $S_2$, $a_2 = \rho a_1-\rho^2$. Suppose $E_2= (\nu_1,-\rho^2+\rho \nu_1) \in S_2$ and $\nu_1 \in [0,2\rho]$. Recall that $E_2=(-2+L\alpha, 1-L \gamma)$. Therefore,
\begin{align}
    \nu_1 &=-2+L\alpha, \label{c1e2x}\\
    -\rho^2+\rho \nu_1 &=1-L \gamma. \label{c1e2y}
\end{align}
Rearranging (\ref{c1e1x}) and (\ref{c1e1y}) yield
\begin{equation}
   \alpha = \frac{2(1-\rho) +\sigma_1}{m} . \label{c1alpha}
\end{equation}
\begin{equation}
    \gamma =\frac{1-\rho^2+\sigma_2}{m}. \label{c1gamma}
\end{equation}
Substituting (\ref{c1alpha}) into (\ref{c1e2x}) leads to
\begin{equation}
    \nu_1 =-2+ \frac{L}{m} \left(2(1-\rho) +\sigma_1 \right). \label{c1nu1}
\end{equation}
Substituting (\ref{c1gamma}) and (\ref{c1nu1}) into (\ref{c1e2y}) and rearranging the terms yields
\begin{equation}
    \kappa =\frac{L}{m}= \frac{(\rho+1)^2}{-3\rho^2+2\rho+1+\rho \sigma_1+\sigma_2}.
\end{equation}
From this equation, it is clear that $\kappa$ achieves its maximum when $\sigma_1 = 0$ and $\sigma_2 = 0$. This indicates that the corresponding value for $\alpha$ is
\begin{equation}
    \alpha = \frac{2(1-\rho)}{m}.
\end{equation}
Also, the requirement that $\nu_1 \in [0,2\rho]$ implies
\begin{align}
    -2+\frac{2(1-\rho)(\rho+1)^2}{-3\rho^2+2\rho+1}& \geq 0,\\
    -2+\frac{2(1-\rho)(\rho+1)^2}{-3\rho^2+2\rho+1}& \leq 2\rho.
\end{align}
However, it follows from these inequalities that
\begin{equation}
    \rho>1,
\end{equation}
which does not satisfy the requirement that $0<\rho<1$. Therefore, Case 1 is excluded from consideration.

\noindent \textbf{Case 2}. Suppose $E_1 = (-2\rho+\sigma_3,\rho^2-\sigma_4)$ where $\sigma_3 \in [0,4\rho)$, $\sigma_4 \in [0,\rho \sigma_3]$ and $E_2 \in S_1 $ as shown in Figure \ref{fig:c2}.
\begin{figure}[!htb] 
\begin{center}
\vspace{-0.4cm}
\includegraphics[scale=0.65]{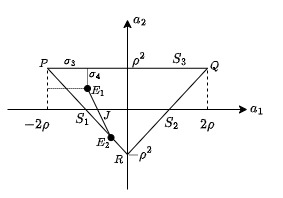}
\end{center} 
\vspace{-0.8cm}
\caption{Illustration of Case 2.}
\vspace{-0.3cm}
\label{fig:c2}
\end{figure}

It follows from (\ref{E1}) that $E_1=(-2+m\alpha,1-m\gamma)$ and hence
\begin{align}
    -2\rho+\sigma_3 =&-2+m\alpha ,\label{c2e1x}\\
    \rho^2-\sigma_4 =& 1-m\gamma .\label{c2e1y}
\end{align}
Also, on the line segment $S_1$, $a_2 = -\rho a_1-\rho^2$. Let $E_2 = (-\nu_2,\rho \nu_2-\rho^2) \in S_1$  for $\nu_2 \in [0,2\rho]$. Recall that $E_2= (-2+ L\alpha,1- L\gamma)$. Therefore,
\begin{align}
   -\nu_2 =& -2+ L\alpha \label{c2e2x},\\
    \rho \nu_2-\rho^2 =& 1- L\gamma \label{c2e2y}.
\end{align}
From rearranging (\ref{c2e1y}), it follows
\begin{equation}
    \gamma = \frac{1-\rho^2+\sigma_4}{m} \label{c2gamma}.
\end{equation}
Substituting (\ref{c2gamma}) into (\ref{c2e2y}) and rearranging, yields
\begin{equation}
    \kappa =\frac{L}{m}= \frac{\rho ^2+1-\rho \nu_2}{1-\rho ^2+\sigma_4}.\label{c2k}
\end{equation}
Since $0<\rho<1$, $1-\rho ^2+\sigma_4 >0$ and $\rho \nu_2>0$,
\begin{equation}
    \kappa \leq \frac{\rho ^2+1 }{1-\rho ^2+\sigma_4}.
\end{equation}
Moreover, $\sigma_4 \in [0,\rho \sigma_3]$. Therefore,
\begin{equation}
    \kappa \leq \frac{\rho ^2+1 }{1-\rho ^2}. \label{c2kineq}
\end{equation}
The upper bound on $\kappa$ in (\ref{c2kineq}) is achieved when $\nu_2 = \sigma_4 = 0$. Hence the maximum value of $\kappa$ in (\ref{c2k}) is
\begin{equation}
    \kappa = \frac{1+\rho ^2}{1-\rho ^2} \label{c2kmax}.
\end{equation}
Substituting $\nu_2 = 0$ and (\ref{c2kmax}) into (\ref{c2e2x}) yields
\begin{equation}
    \alpha = \frac{2(1-\rho^2)}{m(1+\rho^2)}. \label{c2alphamax}
\end{equation}
Also, substituting (\ref{c2alphamax}) into (\ref{c2e1x}) yields
\begin{equation}
    \sigma_3 = 2\rho-2+\frac{2(1-\rho^2)}{1+\rho^2},
\end{equation}
which in turn can be rewritten as 
\begin{equation}
    \sigma_3 = \frac{2\rho(\rho-1)^2}{1+\rho^2}.
\end{equation}
Therefore, it is clear that $\sigma_3>0$ for $0<\rho<1$.

Now, we consider
\begin{align*}
    \sigma_3-4\rho &= \frac{2\rho(\rho-1)^2}{1+\rho^2}-4\rho\\
    &=-\frac{2 \rho  (\rho +1)^2}{\rho ^2+1} <0.
\end{align*}
Hence, $\sigma_3<4\rho$.

Summarising the above cases, we can conclude that the maximum value of $\kappa$ is given by equation (\ref{kmaxp}) with the corresponding values for $\alpha$ and $\gamma$ given by (\ref{kmaxa}) and (\ref{kmaxy}), respectively.  
\end{proof} 

The following lemma pertains to Property 2.
\begin{lemma}
\label{lemspr}
    Given $\kappa = \bar \kappa_{\text{max}}(\rho)$ in (\ref{kmaxp}), $\alpha$ in (\ref{kmaxa}) and $\gamma$ in (\ref{kmaxy}), then the transfer function $H_0(z)$ in (\ref{H0zml}) is strictly positive real.
\end{lemma}
\begin{proof}
Substituting (\ref{kmaxp}), (\ref{kmaxa}), (\ref{kmaxy}) into (\ref{H0zml}) yields
\begin{equation}
    H_0(z) = \frac{z^2-\rho^2}{z^2-\frac{4\rho^2}{1+\rho^2}z+\rho^2}. \label{transferfuncmax}
\end{equation}
To verify that $H_0(z)$ defined in (\ref{transferfuncmax}) is strictly positive real, we make the substitution
\begin{equation}
    z = e^{i \omega} = \cos{\omega} + i\sin{\omega},
\end{equation}
in (\ref{transferfuncmax}), which yields
\begin{align}
    &H_0(e^{i \omega}) \nonumber \\  &=\dfrac{\phi+2i\cos{\omega}\sin{\omega}-2\rho^2}{\phi+2i\cos{\omega}\sin{\omega}-\dfrac{4\rho^2}{1+\rho^2}(\cos{\omega}+i\sin{\omega})} \nonumber\\
    &=\dfrac{\rho+2i\cos{\omega}\sin{\omega}-2\rho^2}{\phi-\dfrac{4\rho^2}{1+\rho^2}\cos{\omega}-i\sin{\omega}(\dfrac{4\rho^2}{1+\rho^2}-2\cos{\omega})} \nonumber\\
    &=\dfrac{(\phi+2i\cos{\omega}\sin{\omega}-2\rho^2)}{ \left( \phi-\dfrac{4\rho^2}{1+\rho^2}\cos{\omega} \right)^2+\left(\sin{\omega}(\dfrac{4\rho^2}{1+\rho^2}-2\cos{\omega}) \right)^2} \nonumber \\ &\quad \quad \times \left(\phi-\dfrac{4\rho^2}{1+\rho^2}\cos{\omega}+i\sin{\omega}(\dfrac{4\rho^2}{1+\rho^2}-2\cos{\omega}) \right) \label{Hzcal}
\end{align}
where $\phi = \cos^2{\omega}-\sin^2{\omega}+\rho^2$.
It is clear that the denominator of (\ref{Hzcal}) is positive. Considering the numerator only
\begin{align*}
    &(\phi+2i\cos{\omega}\sin{\omega}-2\rho^2) \\& \times \left(\phi-\dfrac{4\rho^2}{1+\rho^2}\cos{\omega}+i\sin{\omega}(\dfrac{4\rho^2}{1+\rho^2}-2\cos{\omega}) \right)\\
    &= \left( (\cos{\omega}+i\sin{\omega})^2-\rho^2 \right) \\& \quad \times \left( (\cos{\omega}-i\sin{\omega})^2+\rho^2 -\dfrac{4\rho^2}{1+\rho^2}(\cos{\omega}-i\sin{\omega}) \right)\\
    &= 1-\rho^4+\rho^2(\cos{\omega}+i\sin{\omega})^2-\rho^2(\cos{\omega}-i\sin{\omega})^2 \\& \quad-\left( (\cos{\omega}+i\sin{\omega})^2-\rho^2 \right) \left(\dfrac{4\rho^2}{1+\rho^2}(\cos{\omega}-i\sin{\omega})\right)\\
    &=1-\rho^4+4i\rho^2\cos{\omega}\sin{\omega}\\& \quad -\Big( \dfrac{4\rho^2}{1+\rho^2}(\cos{\omega}+i\sin{\omega}) - \dfrac{4\rho^4}{1+\rho^2}(\cos{\omega}-i\sin{\omega}) \Big)\\
    &= 1-\rho^4+4i\rho^2\cos{\omega}\sin{\omega}\\& \quad -\left(
    \begin{array}{l}
          \dfrac{4\rho^2}{1+\rho^2}\cos{\omega}-\dfrac{4\rho^4}{1+\rho^2}\cos{\omega} \\
          +\dfrac{4\rho^2}{1+\rho^2}i\sin{\omega}+\dfrac{4\rho^4}{1+\rho^2}i\sin{\omega}
    \end{array} \right)\\
    &= 1-\rho^4+4i\rho^2\cos{\omega}\sin{\omega}-\frac{4\rho^2(1-\rho^2)}{1+\rho^2}\cos{\omega} \\& \quad-\frac{4\rho^2(1+\rho^2)}{1+\rho^2}\sin{\omega}\\
    &= 1-\rho^4-\frac{4\rho^2(1-\rho^2)}{1+\rho^2}\cos{\omega}+i4\rho^2\sin{\omega}(\cos{\omega}-1).
\end{align*}
Thus, we need to show that $1-\rho^4-\dfrac{4\rho^2(1-\rho^2)}{1+\rho^2}\cos{\omega} > 0$ for $0<\rho<1$. However
\begin{equation}
    1-\rho^4-\dfrac{4\rho^2(1-\rho^2)}{1+\rho^2}\cos{\omega}> 0,
\end{equation}
if and only if
\begin{equation}
\label{iffdeno}
    \dfrac{4\rho^2(1-\rho^2)}{1+\rho^2}\cos{\omega}<1-\rho^4.
\end{equation}
However, it follows from $0<\rho<1$ that
\begin{equation*}
    \frac{4\rho^2(1-\rho^2)}{1+\rho^2}>0.
\end{equation*}
Hence, (\ref{iffdeno}) is equivalent to
\begin{equation}
    \cos{\omega}<\frac{(1-\rho^4)(1+\rho^2)}{4\rho^2(1-\rho^2)}=\frac{(1+\rho^2)^2}{4\rho^2}. \label{cosomegaprove}
\end{equation}
However, it is straightforward to verifying that
\begin{equation}
    \frac{(1+\rho^2)^2}{4\rho^2} \geq 1 \quad \text{for all} \quad 0<\rho<1
\end{equation}
Given that $-1\leq \cos{\omega}\leq 1$, it follows that (\ref{iffdeno}) must be satisfied.
This completes the proof of the lemma. 
\end{proof}
We are now in a position to prove Theorem \ref{maintheo}.

\textit{Proof of Theorem \ref{maintheo}}:
Given $L>m>0$, $\kappa = \frac{L}{m}$ and $0<\rho<1$, it follows from Lemma \ref{lemkmaxp1} that the maximum value of $\kappa$ such that Property \ref{prop1} is satisfied is given by $\bar \kappa_{\text{max}}(\rho)$ (\ref{kmaxp}), and the corresponding values of $\alpha$ and $\gamma$ are given by (\ref{kmaxa}) and (\ref{kmaxy}) respectively. Moreover, it follows from Lemma \ref{lemspr} that, this $\bar \kappa_{\text{max}}(\rho)$ is also the largest value of $\kappa$ such that Property 2 is satisfied. Thus, the smallest $\rho$ such that Properties \ref{prop1} and \ref{prop2} are satisfied is given by the inverse of the function $\bar \kappa_{\text{max}}(\rho)$. That is 
\begin{equation}
    \rho^*(\kappa) = \sqrt{\frac{\kappa-1}{\kappa+1}},
\end{equation}
with corresponding values of $\alpha$ and $\gamma$ given by
\begin{align}
    \alpha^* &= \frac{2}{m\kappa} = \frac{2}{L},\\
    \gamma^* &=\frac{1}{m(\kappa+1)} = \frac{1}{m+L}.
\end{align}
Also, this implies that $\rho^*(\kappa)$ is the smallest value of $\rho$ such that the characteristic polynomial $\chi(z)$ of $\Tilde{A}_0(\lambda)$ has all its roots inside the disk $|z|<\rho$ and the transfer function $H_0(z)$ is strictly positive real. That is, it follows from Lemmas \ref{rconverge} and \ref{circlecr} that given any algorithm $\Sigma(\alpha,\gamma,m,L)$ which is globally asymptotically convergent via the circle criterion then
\begin{equation*}
    r(\alpha,\gamma,m,L) \geq \rho^*(\kappa).
\end{equation*}
This completes the proof. $\hfill \blacksquare$

\section{Illustrative Example}\label{sec4}
This section provides an illustrative example of the use of algorithm (\ref{uprulebeta1}) in a time of arrival (TOA) based source localization problem. We consider a TOA based localization problem similar to the problem considered in \cite{MandyTOA}. Let $x^*(t) \in \mathbb{R}^n$ denote the real source location at time $t$. Consider an array of $m$ sensors and let $s_i\in \mathbb{R}^n$ be the location of the $i$th sensor. Let $r_i$ denote the noisy range measurement between the source and the $i$th sensor
\begin{equation}
    r_i(t) = \|x^*(t)-s_i\|+ \epsilon_i(t), \quad i = 1,\dots,m ,
\end{equation}
where $\epsilon_i(t)$ is the measurement noise. The source location is estimated by the minimization of the least square cost function
\begin{equation}
\label{LQTOA}
    f(x,t):=  \sum_{i=1}^{m} (\|x(t)-s_i\|-r_i(t))^2
\end{equation}
Here, $x(t) = [x_1(t), x_2(t)]^T $ denotes the estimate of the true source location at time step $t$. Also, the gradient and Hessian of the cost function defined in (\ref{LQTOA}) are calculated to be
\begin{equation*}
    \nabla_x f(x(t),t) = \sum_{i=1}^{m} 2\left( x(t)-s_i-r_i(t) \right)\frac{x(t)-s_i}{\|x(t)-s_i\|},
\end{equation*}
\begin{align*}
    \nabla_x^2 f(x(t),t) &= \sum_{i=1}^{m} 2 \big(I - r_i (\frac{I}{\|x(t)-s_i\|}\nonumber \\&\quad-\frac{(x(t)-s_i)(x(t)-s_i)^T}{\|x(t)-s_i\|^3}) \big).
\end{align*}

We now assume that the true source location can be perfectly measured without noise, which yields
\begin{equation}
    \epsilon_i(t) = 0.
\end{equation}
Therefore, we can rewrite (\ref{LQTOA}) as
\begin{equation}
\label{LQTOA1}
    f(x,t):=  \sum_{i=1}^{m} (\|x(t)-s_i\|-\|x^*(t)-s_i\|)^2.
\end{equation}
The Hessian of equation (\ref{LQTOA1}) does not satisfy the condition (\ref{classfuncmL}) globally and uniformly in time. Hence, we only consider a specific time horizon of interest.

Assume that there are three sensors ($m=3$) whose coordinates are $s_1 = [1, 0.8]^T$, $s_2 = [1, -1]^T$ and $s_3 = [0, -0.5]^T$ with initial true source location $x^*(0)=[-9,10]^T$. In order to apply our algorithm, we consider the situation where the source moves linearly at a constant speed and direction. This yields 
\begin{equation}
\label{TOAspeed}
    x^*(t) = x^*(0) +at,
\end{equation}
where $a = [0.01,-0.01]^T$ represents the velocity vector of the source.

\begin{figure}[!htb] 
    \vspace{-0.3cm}
    \centering
    \includegraphics[width=\linewidth]{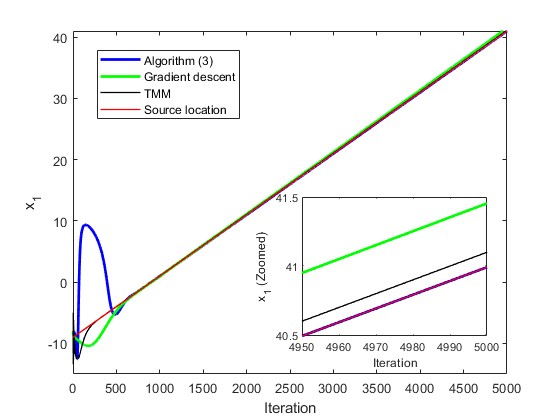}
    \vspace{-0.8cm}
    \caption{Comparison between the algorithm (\ref{uprulebeta1}), the gradient descent algorithm \cite{Boyd_Vandenberghe_2004} and the triple momentum method \cite{7967721} applied to the TOA-based localization problem (\ref{LQTOA}), with $L = 6$ and $m = 0.1$, to estimate the true source location $x_1$ at different iteration steps $t$.}
    \label{fig:x1}
    \vspace{-0.2cm}
\end{figure}

\begin{figure}[!htb] 
    \vspace{-0.2cm}
    \centering
    \includegraphics[width=\linewidth]{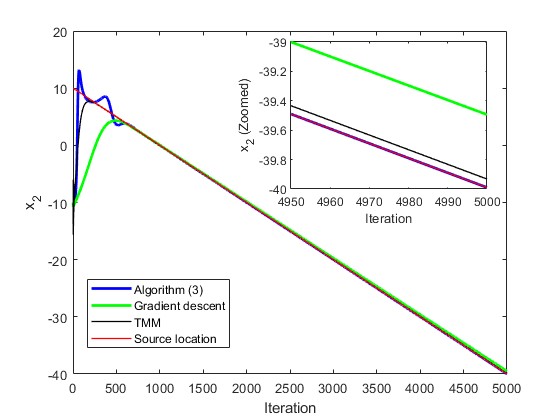}
    \vspace{-0.8cm}
    \caption{Comparison between the algorithm (\ref{uprulebeta1}), the gradient descent algorithm \cite{Boyd_Vandenberghe_2004} and the triple momentum method \cite{7967721} applied to the TOA-based localization problem (\ref{LQTOA}), with $L = 6$ and $m = 0.1$, to estimate the true source location $x_2$ at different iteration steps $t$.}
    \label{fig:x2}

\end{figure}

Figures \ref{fig:x1} and \ref{fig:x2} compare the true position of the source and the estimated position obtained using algorithm (\ref{uprulebeta1}), the gradient descent algorithm with the step size of $\frac{2}{L+m}$ and the triple momentum method (TMM) \cite{7967721}, respectively, with the estimation initialized at $x(0) = [-8,-10]^T$. The red line demonstrates the true source location, which moves linearly corresponding to equation (\ref{TOAspeed}). The blue line represents the source location estimated using Algorithm (\ref{uprulebeta1}). Although a significant transient can be observed initially, the algorithm ultimately succeeds in estimating the true source location with zero steady-state error. This behavior demonstrates that the algorithm converges with zero steady-state error, indicating its efficacy in reaching an accurate solution despite initial fluctuations. The green line and blue line in Figures \ref{fig:x1} and \ref{fig:x2} illustrate the performance of the gradient descent algorithm and TMM in estimating the true source location. 

\begin{figure}[!htb] 
    \vspace{-0.3cm}
    \centering
    \includegraphics[scale=0.48]{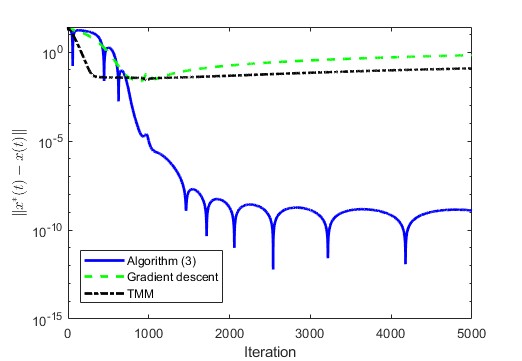}
    \vspace{-0.8cm}
    \caption{Difference between the true source location and its estimate using algorithm (\ref{uprulebeta1}), the gradient descent algorithm \cite{Boyd_Vandenberghe_2004} and the triple momentum method \cite{7967721}.}
    \label{fig:error}
    \vspace{-0.6cm}
\end{figure}

In Figure \ref{fig:error}, both gradient descent algorithm and TMM demonstrate a non-zero steady-state error in estimating the true source location. This discrepancy may be attributed to the fact that the gradient descent algorithm and the TMM algorithm are not designed for estimating parameters in time-varying scenarios and only contain single integration in the algorithm, leading to steady-state errors in convergence to the true values. Additionally, TMM is designed for optimizing strongly convex functions, which does not align with the characteristics of problem (\ref{LQTOA1}). Figure \ref{fig:error} also verifies that the estimation using algorithm (\ref{uprulebeta1}) achieves zero steady-state error, whereas the gradient descent method and the TMM do not.

\section{CONCLUSIONS}\label{sec5}
In this paper, we proposed an online optimization algorithm for solving a class of nonconvex optimization problems with sector bounded gradients. Our results demonstrate that the proposed algorithm  converges globally to a linearly varying optimal point with zero steady-state error. This property was verified through an illustrative example involving a TOA based localization problem.








\bibliographystyle{IEEEtran}
\bibliography{IEEEabrv,irpnew}

\end{document}